\documentclass{amsart}

\usepackage{amssymb}
\usepackage{amsmath}
\usepackage{amsthm}

\usepackage{graphicx}

\usepackage{amscd}
\usepackage[cmtip,all]{xy}

\newtheorem{theorem}{Theorem}[section]
\newtheorem{lemma}[theorem]{Lemma}

\theoremstyle{remark}
\newtheorem{remark}[theorem]{Remark}

\begin{document}
\title[Ring $K(s)^*(BG)$, $G=G_{38},\cdots,G_{41}$]{Morava $K$-theory rings for the groups $G_{38},\cdots,G_{41}$ of order 32 }
\author{Malkhaz Bakuradze and Mamuka Jibladze}
\address{Iv. Javakhishvili Tbilisi State University, Faculty of Exact and Natural Sciences}
\email{malkhaz.bakuradze@tsu.ge,  jib@rmi.ge}
\thanks{First author was supported by Volkswagen Foundation, Ref. 1/84 328 and Rustaveli Foundation grant DI/16/5-103/12 . Second author was supported by the STCU grant 5622 and Rustaveli Foundation grant 09/23.}


\date{}


\begin{abstract}
B. Schuster \cite{SCH1} proved that the $mod$ 2 Morava $K$-theory $K(s)^*(BG)$ is evenly generated for all groups $G$ of order 32. For the four groups $G$  with the numbers 38, 39, 40 and 41 in the Hall-Senior list \cite{H}, the ring $K(2)^*(BG)$ has been shown to be generated as a $K(2)^*$-module by transferred Euler classes.  In this paper, we show this for arbitrary $s$ and compute the ring structure of $K(s)^*(BG)$. Namely, we show that $K(s)^*(BG)$ is the quotient of a polynomial ring in 6 variables over $K(s)^*(pt)$ by an ideal for which we list explicit generators.
\end{abstract}

\maketitle{}


{\centering\section{Introduction and Statements}}

Let $K(s)^*$, $s>1$, be the $s$-th Morava $K$-theory at 2.
In this paper we compute the ring structure of $K(s)^*(BG)$ for the four groups $G=G_{38},...,G_{41}$ from the Hall-Senior list \cite{H}, by showing that $K(s)^*(BG)$ is the quotient of a polynomial ring $K(s)^*(pt)[a,b,c,x_2,y_2,T]$ by a certain ideal $R$ for which we give explicit generators.

A finite group $G$ is said to be good \cite{HKR} if $K(s)^*(BG)$ is generated as a $K(s)^*$-module by transfers of Euler classes of complex representations.
Special effort was needed to find an example of a group not good in this sense \cite{K}. For the additive structure, the principal calculational tool is the Atiyah-Hirzebruch spectral sequence \cite{BU1,BU2} and the Serre SS \cite{K}. Even if the additive structure is calculated, the multiplicative structure is still a delicate task. It is not always determined by representation theory, i.e., $G$ does not have exact Chern approximation in the terminology of Strickland \cite{ST}.  Also the presentation of $K(s)^*(BG)$ in terms of the formal group law and splitting principle \cite{KA} is not always convenient. This clearly indicates that the part of the relations which can be derived from the properties of the transfer should play decisive role in determining the whole ring structure.

In the current paper we will consider four groups  $G=G_{38},\cdots,G_{41}$ of order 32 from the Hall-Senior list \cite{H}. It is proved in \cite{SCH1} that $K(s)^*(BG)$ is evenly generated and for $s=2$ is generated by Euler classes and transferred Euler classes. One consequence of our main theorem below is that this is true for any $s$.
We obtain generators for the ideal $R$ above by using the formula for transferred Euler class from \cite{BP} and follow a certain plan, which proved to be sufficient to handle the 2-groups $D$,\,$SD$,\,$QD$,\,$Q$ \cite{BV}, \cite{B2} and modular $p$-groups \cite{B1}.
For a discussion of the ring  structure of all other groups of order 32 see \cite{SCH1}, \cite{SCH2}.
\medskip

Let $G$ be one of the groups
\begin{multline*}
\shoveleft{}\\
\shoveleft{G_{38}=\langle \mathbf{a},\mathbf{b},\mathbf{c} \mid  \mathbf{a}^4=\mathbf{b}^2=\mathbf{c}^4=[\mathbf{a},\mathbf{b}]=1, \mathbf{c}\mathbf{a}\mathbf{c}^{-1}=\mathbf{a}\mathbf{c}^2, \mathbf{c}\mathbf{b}\mathbf{c}^{-1}=\mathbf{a}^2\mathbf{b}\rangle ,}\\
\shoveleft{G_{39}=\langle \mathbf{a},\mathbf{b},\mathbf{c} \mid  \mathbf{a}^4=\mathbf{b}^4=\mathbf{c}^2=[\mathbf{a},\mathbf{b}]=1,
\mathbf{c}\mathbf{a}\mathbf{c}=\mathbf{a}^3, \mathbf{c}\mathbf{b}\mathbf{c}=\mathbf{a}^2\mathbf{b}^3 \rangle ,}\\
\shoveleft{G_{40}=\langle \mathbf{a},\mathbf{b},\mathbf{c} \mid  \mathbf{a}^4=\mathbf{b}^4=1, \mathbf{c}^2=\mathbf{b}^2, [\mathbf{a},\mathbf{b}]=1, \mathbf{c}^{-1}\mathbf{a}\mathbf{c}=\mathbf{a}^3, \mathbf{c}^{-1}\mathbf{b}\mathbf{c}=\mathbf{a}^2\mathbf{b}^3 \rangle ,}\\
\shoveleft{G_{41}=\langle \mathbf{a},\mathbf{b},\mathbf{c} \mid  \mathbf{a}^4=\mathbf{b}^4=\mathbf{c}^2=[\mathbf{a},\mathbf{b}]=1,  \mathbf{c}\mathbf{a}\mathbf{c}=\mathbf{a}^3\mathbf{b}^2, \mathbf{c}\mathbf{b}\mathbf{c}=\mathbf{a}^2\mathbf{b} \rangle .}\\
\end{multline*}
Let $H$ be the maximal abelian subgroup of index two $\langle \mathbf{a}, \mathbf{b}, \mathbf{c}^2\rangle \cong C_4 \times C_2 \times C_2$ for $G=G_{38}$ and $\langle \mathbf{a}, \mathbf{b}\rangle  \cong C_4 \times C_4$ for all other cases.
Let $\lambda$, $\mu$ and $\nu$ denote complex line bundles over $BH$. For $H \lhd G_{38}$, let
$$\lambda(\mathbf{a})=i, \mu(\mathbf{b})=\nu(\mathbf{c}^2)=-1, \lambda(\mathbf{b})=\lambda(\mathbf{c}^2)=\mu(\mathbf{a})=\mu(\mathbf{c}^2)=\nu(\mathbf{a})=\nu(\mathbf{b})=1,$$
be the pullbacks of the canonical complex line bundles along the projections onto the first, second and third factor of $H$ respectively. For all other cases, let

$$\lambda(\mathbf{a})=\nu(\mathbf{b})=i, \lambda(\mathbf{b})=\nu(\mathbf{a})=1,$$
be the pullbacks of the canonical complex line bundles along the projections onto the first and second factor of $H$ respectively.

The quotient of $G$ by the center is isomorphic to $ C_2 \times C_2 \times C_2$. The projections on the three factors induce three line bundles $\alpha$, $\beta$ and $\gamma$ respectively.

Let us denote Chern classes by

\begin{equation*}
x_i=
\begin{cases}
c_i(Ind_H^G(\lambda))& \text{for $G=G_{39}, G_{40}$}\\
c_i(Ind_H^G(\nu)) & \text{for $G=G_{38}, G_{41}$};
\end{cases}
\end{equation*}

\medskip

\begin{equation*}
y_i=
\begin{cases}
c_i(Ind_H^G(\nu)) & \text{for $G=G_{39}, G_{40}$}\\
c_i(Ind_H^G)(\lambda) & \text{for $G=G_{38},G_{41}$};
\end{cases}
\end{equation*}

\medskip

\begin{equation*}
a=\begin{cases}
c_1(\alpha) & \text{for $G=G_{38}, G_{41}$}\\
c_1(\alpha \beta) & \text{for $G=G_{39}, G_{40}$};
\end{cases}
\end{equation*}

\medskip

\begin{equation*}
b=\begin{cases}
c_1(\beta) & \text{for $G=G_{38}, G_{39}, G_{40}$}\\
c_1(\alpha \beta) & \text{for $G=G_{41}$};
\end{cases}
\end{equation*}

\medskip

\begin{equation*}
c=c_1(\gamma), \,\, \text{for all cases.}
\end{equation*}

\medskip

Let $Tr^*: K(s)^*(BH)\to K(s)^*(BG)$ be the transfer homomorphism \cite{A} associated to the double covering $\rho: BH \to BG$ and let
$$
T=Tr^*(c_1(\lambda)c_1(\nu)).
$$

\medskip

 Note that by \cite{JW}, $K(s)^*(pt)$ is the Laurent polynomial ring in one variable, which is usually denoted in our situation by $\mathbb{F}_2[v_s,v_s^{-1}]$, where $\mathbb{F}_2$ is the field of 2 elements.

Our main result is the following

\bigskip

\medskip

\begin{theorem}
\label{thm:G}  Let $G$ be one of the groups $G_{38},...,G_{41}$. Then

i) $K(s)^*(BG)\cong K(s)^*[a,b,c,x_2,y_2,T]/R$, where the ideal $R$ is generated by

$a^{2^s}$, $b^{2^s}$, $c^{2^s}$,

$c(c+x_1+v_s\sum_{i=1}^{s-1}c^{2^s-2^i}x_2^{2^{i-1}})$,
\,\,$c(c+y_1+v_s\sum_{i=1}^{s-1}c^{2^s-2^i}y_2^{2^{i-1}})$,

$a(a+x_1+v_s\sum_{i=1}^{s-1}a^{2^s-2^i}x_2^{2^{i-1}})$,
\,\,$b(b+y_1+v_s\sum_{i=1}^{s-1}b^{2^s-2^i}y_2^{2^{i-1}})$,

$(c+x_1+v_s\sum_{i=1}^{s-1}c^{2^s-2^i}x_2^{2^{i-1}})(b+y_1+v_s\sum_{i=1}^{s-1}b^{2^s-2^i}y_2^{2^{i-1}})+v_sb^{2^s-1}T$,

$(c+y_1+v_s\sum_{i=1}^{s-1}c^{2^s-2^i}y_2^{2^{i-1}})(a+x_1+v_s\sum_{i=1}^{s-1}a^{2^s-2^i}x_2^{2^{i-1}})+v_sa^{2^s-1}T$,

$T^2+Tx_1y_1+x_2y_1(c+y_1+v_s\sum_{i=1}^{s-1}c^{2^s-2^i}y_2^{2^{i-1}})+x_1y_2(c+x_1+v_s\sum_{i=1}^{s-1}c^{2^s-2^i}x_2^{2^{i-1}})$,

$T(a+x_1+v_s\sum_{i=1}^{s-1}a^{2^s-2^i}x_2^{2^{i-1}})+v_sa^{2^s-1}x_2(c+y_1)$,

$T(b+y_1+v_s\sum_{i=1}^{s-1}b^{2^s-2^i}y_2^{2^{i-1}})+v_sb^{2^s-1}y_2(c+x_1)$,\,\,\,$cT$, and

\begin{equation*}
v_s^2x_2^{2^s}+
\begin{cases}
a^2+b^2+ac+v_sabc^{2^s-1}& \text {for $G=G_{39},G_{40},G_{41}$}\\
c^2+ac & \text{\,\,\, $G=G_{38}$}
\end{cases}
\end{equation*}

\begin{equation*}
v_s^2y_2^{2^s}+
\begin{cases}
a^2+bc+v_sabc^{2^s-1}& \text{for $G=G_{38},G_{41}$}\\
b^2+bc & \text{\,\,\,$G=G_{39}$} \\
b^2+c^2+bc & \text{\,\,\,$G=G_{40}$},
\end{cases}
\end{equation*}

where
\begin{equation*}
x_1=v_s(x_2+v_sx_1x_2^{2^{s-1}})^{2^{s-1}}+
\begin{cases}
a & \text{for $G=G_{38}$}\\
b+c+v_s(bc)^{2^{s-1}} & \text{\,\,\,$G=G_{39},G_{40},G_{41};$}
\end{cases}
\end{equation*}

\begin{equation*}
y_1=v_s(y_2+v_sy_1y_2^{2^{s-1}})^{2^{s-1}}+
\begin{cases}
c& \text{for $G=G_{39}$} \\
0& \text{\,\,\, $G=G_{40}$}\\
a+b+c+ & \\
v_s(ab+bc+ac)^{2^{s-1}}& \text{$G=G_{38},G_{41}.$}
\end{cases}
\end{equation*}

ii) Some other relations are

$a^2c=ac^2,\,\,\,$ $b^2c=bc^2,\,\,\,$ $x_1^{2^{s}}=a^{2^{s-1}}c^{2^{s-1}},\,\,\,$ $y_1^{2^{s}}=b^{2^{s-1}}c^{2^{s-1}}.$
\end{theorem}

\bigskip

\bigskip

The rest of the paper is organized as follows. Section 2 presents some preliminaries. In Section 3 we treat the representation theory of the groups under consideration. In Section 4 we derive the relations of Theorem \ref{thm:G}. Section 5 is devoted to the most difficult part of the proof of Theorem \ref{thm:G}. Namely for each of our groups (see Lemma \ref{lem:free4module}, Lemma \ref{basis-39-40} for $G_{39}$, $G_{40}$ and Lemma \ref{lem:H38free4module}, Lemma \ref{basis38-41} for $G_{38}$, $G_{41}$) we prove that certain monomials in $a,b,c, x_1,x_2,y_1,y_2, T $ form a basis of $K(s)^*(BG)$ as a free $K(s)^*$-module. It follows $c,a,b,x_2,y_2,T$ are $K(s)^*$-algebra generators as $x_1$ and $y_1$ are decomposable in these elements. Finally we prove that the relations in Section 4 provide a complete set of defining relations. For the reader's convenience Section 6 discusses some papers on the subject.

\section{Preliminaries}

Let $H\lhd G$ be of index 2. Consider the double covering $\pi : BH \to BG$. Let
$$Tr^*_{\pi}=Tr^*(H,G): K(s)^*(BH)\to K(s)^*(BG)$$
be the associated transfer homomorphism induced by the stable transfer map \cite{A}, \cite{KP}, \cite{D}.
We will need the following transfer formula from \cite{BP}.

\medskip

Let $\xi \to BH$ be a complex line bundle and $\xi_{\pi}=Ind_H^G(\xi)$ be its Atiyah transfer. Then
\begin{equation}
\label{eq:tr}
c_1(\xi_{\pi})=
c_1(\psi)+v_s\sum_{i=1}^{s-1}c_1(\psi)^{2^s-2^i}c_2(\xi_{\pi})^{2^{i-1}}
+Tr^*_{\pi}(c_1(\xi)),
\end{equation}
where $\psi \rightarrow BG$ is the pullback of the canonical line bundle over $B\mathbb{Z}/2$ along the map $BG \rightarrow B\mathbb{Z}/2$ classifying $\pi$.

\medskip

Let

\begin{equation*}
u=
\begin{cases}
c_1(\lambda)& \text{for $G=G_{39}, G_{40}$}\\
c_1(\nu) & \text{for $G=G_{38}, G_{41}$}
\end{cases}
\end{equation*}
and
\begin{equation*}
v=
\begin{cases}
c_1(\lambda)& \text{for $G=G_{38}, G_{41}$}\\
c_1(\nu) & \text{for $G=G_{39}, G_{40}$}.
\end{cases}
\end{equation*}

For the Chern classes $u$, $v$ and $Tr^*=Tr^*(H,G)$ the formula \eqref{eq:tr} implies

\begin{equation}
\label{Tru}
Tr^*(u)=c+x_1+v_s\sum_{i=1}^{s-1}c^{2^s-2^i}x_2^{2^{i-1}}
\end{equation}
and

\begin{equation}
\label{Trv}
Tr^*(v)=c+y_1+v_s\sum_{i=1}^{s-1}c^{2^s-2^i}y_2^{2^{i-1}}.
\end{equation}

\medskip

One has the following approximation formula for the formal group law in Morava $K$-theory (\cite{BV}, Lemma 2.2 ii)).

\begin{equation}
\label{eq:FGL}
F(x,y)=x+y+v_s\Phi(v_s,x,y)^{2^{s-1}},
\end{equation}
where $\Phi(v_s,x,y)=xy+v_s(xy)^{2^{s-1}}(x+y)\,\,\, modulo \,\, (xy)^{2^{s-1}}(x+y)^{2^{s-1}}.$

\bigskip

We also will need the following

\begin{lemma}
\label{lem:zeta^2}
The tensor square of a complex plane vector bundle $\zeta$ has the following total Chern class
$C(\zeta^{\otimes 2})=(1+c_1^2(det\zeta))(1+v_sc_1^{2^s}(\zeta)+v_s^2c_2^{2^s}(\zeta)).$
\end{lemma}
\begin{proof}
Use the splitting principle and write formally
$$\zeta=\xi_1+\xi_2$$ and
$$c_1(\zeta)=t_1+t_2; \,\,\, c_2(\zeta)=t_1t_2.$$
We have that the $i$-th Chern class ($i=1,2,3,4$) on the right hand side of the bundle relation
\begin{equation}
(\xi_1+\xi_2)\otimes(\xi_1+\xi_2)=\xi_1^2+\xi_2^2+2\xi_1\otimes \xi_2
\end{equation}
is the $i$-th elementary symmetric function in $F(t_1,t_1),F(t_2,t_2),F(t_1,t_2),F(t_1,t_2)$.
That is
$$c_i(\zeta^2)=\sigma_i(F(t_1,t_1),F(t_2,t_2),F(t_1,t_2),F(t_1,t_2)).$$
Hence we have for the first Chern class
$$c_1(\zeta^2)=v_st_1^{2^s}+v_st_2^{2^s}=v_sc_1^{2^s}(\zeta).$$
For the second Chern class we have
$$c_2(\zeta^2)=v_s t_1^{2^s}v_st_2^{2^s}+F^2(t_1,t_2)=v_s^2c_2^{2^s}(\zeta)+c_1^2(det\zeta).$$
Similarly, the third and fourth Chern classes are
$$v_s(t_1^{2^s}+t_2^{2^s})F^2(t_1,t_2)=v_sc_1^{2^s}(\zeta)c_1^2(det\zeta)$$
and
$$v_st_1^{2^s}v_st_2^{2^s}F^2(t_1,t_2)=v_s^2c_2^{2^s}(\zeta)c_1^2(det\zeta)$$
respectively.
\end{proof}

\section{Bundle relations}

Let us give some relations of bundles over $BG$ we will need. We omit the proofs since they are completely standard and easily follow from the definitions and from Frobenius reciprocity of the transfer in complex $K$-theory.

Let $\rho:BH \rightarrow BG$ be the double covering $\rho=\rho(H,G)$ and let $\lambda_{!}=Ind^G_H(\lambda)$ and $\nu_{!}=Ind^G_H(\nu)$ in each of the four cases.

\bigskip

\noindent{$\bf G=G_{38}$; $H=\langle \mathbf{a},\mathbf{b},\mathbf{c}^2\rangle$}.

\medskip

Determinants.

\noindent{\,}$det(\lambda_!)=\alpha\beta\gamma, \,\, det(\nu_!)=\alpha$ and

\medskip

Restrictions.

\noindent{i)} $\rho^*\alpha=\lambda^ {2},$  $\rho^*\beta=\mu,$  $\rho^*\gamma=1,$

\medskip

\noindent{ii)} $\rho^*\lambda_!=\lambda +\lambda \mu ,$ $\rho^*\nu_!=\nu + \lambda^{2} \nu ;$

\medskip

Product relations.

\noindent{iii)} $\beta \lambda_!=\lambda_!,$ $\gamma \lambda_!=\lambda_!$

\medskip

\noindent{iv)}  $\alpha \nu_!=\nu_!,$ $\gamma \nu_!=\nu_!;$

\medskip

\noindent{v)} $(\nu_{!})^2=1 + \alpha + \gamma + \alpha \gamma.$

\medskip

\noindent{vi)} $(\lambda_{!})^2=\alpha + \alpha \gamma + \alpha \beta + \alpha \beta \gamma.$

\bigskip

The first relation of iv) suggests that $\nu_!$ should be also transferred from some line bundle for the $2$-covering corresponding to $\alpha.$ Namely we have

\medskip

\begin{lemma}
\label{lem:NU38}
Let $A=\langle \mathbf{a}^2,\mathbf{b},\mathbf{c}\rangle $ and let $\nu'$ be represented by $\nu'(\mathbf{a}^2)=1,$ $\nu'(\mathbf{b})=1,$ $\nu'(\mathbf{c})=i$.
Then $Ind_H^G(\nu)=Ind_{A}^G(\nu').$
\end{lemma}

Similarly the second relation of iii) suggests

\begin{lemma}
\label{lem:LAMBDA38}
Let $B=\langle \mathbf{a},\mathbf{c}\rangle \cong \langle \mathbf{c} \rangle \rtimes \langle \mathbf{a} \rangle$ and let $\lambda'$ be represented by $\lambda'(\mathbf{a})=i,$ $\lambda'(\mathbf{c})=1$. Then $Ind_H^G(\lambda)=Ind_B^G(\lambda').$
\end{lemma}

\bigskip

\bigskip

\noindent{$\bf G=G_{39}$, $H=\langle \mathbf{a}, \mathbf{b}\rangle\cong C_4 \times C_4$.}

\medskip

Determinants.

\noindent{} $det(\lambda_!)=\beta \gamma, \,\, det(\nu_!)=\gamma$.

\medskip

Restrictions.

\noindent{i)} $\rho^*\alpha=\lambda^{2},$  $\rho^*\beta=\nu^2,$  $\rho^*\gamma=1;$

\medskip

\noindent {ii)} $\rho^* \lambda_!=\lambda +\lambda^3 \nu^2 $, $\rho^*\nu_!=\nu + \nu^3.$

\medskip

Product relations.

\noindent{iii)} $\alpha \beta \lambda_!=\lambda_! ,$  $\gamma \lambda_!=\lambda_!;$


\medskip

\noindent{iv)} $\beta \nu_!=\nu_!,$  $ \gamma \nu_!=\nu_!;$

\medskip

\noindent{v)} $(\lambda_{!})^2=\alpha + \beta + \alpha \gamma + \beta \gamma.$

\medskip

\noindent{vi)} $(\nu_{!})^2=1+\beta+\gamma + \beta \gamma.$

\bigskip

The first relation of iv) suggests that $\nu_!$ should be also transferred from some line bundle for the $2$-covering corresponding $\beta.$ Namely we have

\medskip

\begin{lemma}
\label{lem:NU39}
Let $B=\langle \mathbf{a},\mathbf{b}^2,\mathbf{c}\rangle $ and let $\nu'$ be represented by $\nu'(\mathbf{a})=1,$ $\nu'(\mathbf{b}^2)=-1, $ $\nu'(\mathbf{c})=1,$
Then $Ind_H^G(\nu)=Ind_{B}^G(\nu').$
\end{lemma}





\medskip




The second relation of iii) suggests

\begin{lemma}
\label{lem:LAMBDA39}
Let $AB=\langle \mathbf{a}^2,\mathbf{a}\mathbf{b},\mathbf{c}\rangle$ and let $\lambda'$ be represented by $\lambda'(\mathbf{a}^2)=-1,$ $\lambda'(\mathbf{a}\mathbf{b})=i,$ $\lambda'(\mathbf{c})=1$.
Then $Ind_H^G(\lambda)=Ind_{AB}^G(\lambda').$
\end{lemma}








\bigskip

\noindent{$\bf G=G_{40}$,\, $H=\langle \mathbf{a}, \mathbf{b}\rangle\cong C_4 \times C_4$.

\medskip

$G_{40}$ has the same character table as $G_{39}$. The only difference is in the determinants of $\lambda_{!}$ and $\nu_{!}$, while restrictions, products, and the two lemmas above are the same.

$$det(\lambda_!)=\beta \gamma, \,\, det(\nu_!)=1.$$

\bigskip

\noindent{$\bf G=G_{41}$, \, $H=\langle \mathbf{a}, \mathbf{b}\rangle\cong C_4 \times C_4$.}

\medskip

Determinants.

\noindent{} $det(\lambda_!)=\beta \gamma, \,\, det(\nu_!)=\alpha \beta \gamma.$

\medskip

Restrictions.

\noindent{i)} $\rho^*\alpha=\lambda^ {2},$  $\rho^*\beta=\nu^2,$  $\rho^*\gamma=1;$

\medskip

\noindent{ii)} $\rho^* \lambda_!=\lambda +\lambda^3 \nu^2 $, $\rho^*\nu_!=\nu + \lambda^2 \nu.$

\medskip

Product relations.

\noindent{iii)} $\gamma \lambda_!=\alpha \beta \lambda_!=\lambda_!;$


\medskip

\noindent{iv)} $\alpha \nu_!=\gamma \nu_!=\nu_!$.

\medskip

\noindent{v)} $(\lambda_{!})^2=\alpha + \beta + \alpha \gamma + \beta \gamma.$

\medskip

\noindent{vi)} $(\nu_{!})^2=\beta + \alpha \beta+ \beta \gamma + \alpha \beta \gamma.$

\bigskip

\begin{lemma}
\label{lem:LAMBDA41}
We can replace the group $G_{40}$ by $G_{41}$ in Lemma \ref{lem:LAMBDA39}.
\end{lemma}


Also for $\nu_!$ one has



\begin{lemma}
\label{lem:NU41}
$\nu_!=Ind^G_A(\nu'),$  where $A={\langle \mathbf{a}^2,\mathbf{b},\mathbf{c} \rangle},$
and $\nu'(\mathbf{a}^2)=1,$ $\nu'(\mathbf{b})=i,$ $\nu'(\mathbf{c})=1.$
\end{lemma}

\bigskip

\bigskip

\section{Relations of Theorem \ref{thm:G}}

Clearly the relations
$$
a^{2^s}=b^{2^s}=c^{2^s}=0.
$$
are immediate consequences of the bundle relations $\alpha^2=\beta^2=\gamma^2=1$ for all cases.

The 4th and 5th relations follow from \eqref{Tru} and \eqref{Trv} respectively.

For the 6th relation
\begin{equation*}
a(a+y_1+v_s\sum_{i=1}^{s-1}a^{2^s-2^i}y_2^{2^{i-1}})=0
\end{equation*}
consider the double covering $\rho:BH \to BG$ in each of the four cases and apply formula \eqref{eq:tr} and Lemma \ref{lem:NU38}, \ref{lem:NU39}, \ref{lem:NU39} or \ref{lem:LAMBDA41} for $G_{38}$, $G_{39}$, $G_{40}$, or $G_{41}$ respectively. For example if $G=G_{38}$ formula \eqref{eq:tr} and Lemma \ref{lem:NU38} imply that the second factor of the relation is the transfer of $c_1(\nu')$
$$Tr^*(c_1(\nu'))=(a+y_1+v_s\sum_{i=1}^{s-1}a^{2^s-2^i}y_2^{2^{i-1}}).$$
Then
$aTr^*(c_1(\nu'))=Tr^*(\rho^*(a)c_1(\nu'))=Tr^*(0\cdot c_1(\nu'))=0$.
\qed

\medskip

Similarly for the 7th relation
\begin{equation*}
b(b+x_1+v_s\sum_{i=1}^{s-1}b^{2^s-2^i}x_2^{2^{i-1}})=0
\end{equation*}
apply formula \eqref{eq:tr} and Lemma \ref{lem:LAMBDA38}, \ref{lem:LAMBDA39}, \ref{lem:LAMBDA39} or \ref{lem:NU41} for $G_{38}$, $G_{39}$, $G_{40}$, or $G_{41}$ respectively.
\qed

\medskip

Now note that the 4th and 6th relations imply $a^2c=ac^2$, the first relations of Theorem \ref{thm:G} ii). For this multiply the 4th relation by $a$ and the 6th relation by $c$. The sum of these terms equals $a^2c+ac^2$ up to an invertible factor.  Similarly the 5th and 7th relation imply $b^2c=bc^2$, the second relation of Theorem \ref{thm:G} ii).

\medskip

For the decompositions of $v_s^2x_2^{2^s}$, $v_s^2y_2^{2^s}$, (also for the formulas for $x_1^{2^s}$ and $y_1^{2^s}$ of Theorem \ref{thm:G} ii) we need the material of Section 3. Namely we have to apply Lemma \ref{lem:zeta^2} to all induced representations given in Section 3 and take into account that their determinants can written in terms of the bundles $\alpha, \beta, \gamma$.
For example for $G=G_{38}$

\begin{equation*}
v_s^2x_2^{2^{s}}=c^2+ac, \,\,\, x_1^{2^s}=a^{2^{s-1}}c^{2^{s-1}}
\end{equation*}
and
\begin{equation*}
v_s^2y_2^{2^s}=a^2+bc+v_sabc^{2^s-1}, \,\,\, y_1^{2^s}=b^{2^{s-1}}c^{2^{s-1}}
\end{equation*}
are the consequences of product relations v) and vi).
Let us prove first two relations. Equate Chern classes in the bundle relation of v). Then for the first Chern classes we get
$$
v_sx_1^{2^s}=a+c+a+c+v_sa^{2^{s-1}}c^{2^{s-1}}=v_sa^{2^{s-1}}c^{2^{s-1}}.
$$

For the decomposition of $v_s^2x_2^{2^s}$ apply the equation for the second Chern classes:

\begin{align*}
v_s^{2}x_2^{2^s}&=c_2(\nu_{!}^2)+c_1(det\,\nu_{!})^2 && \\
&=c_2(1+\alpha+\gamma+\alpha \gamma)+c_1(\alpha)^2=ac+(a+c)F(a,c)+a^2 && \\
&=c^2+ac+v_s(a+c)(ac)^{2^{s-1}} && \\
&=c^2+ac &&\\
\end{align*}

since $a^2c=ac^2$.

\medskip

Note also

\begin{equation}
\label{a^kc^l}
a^kc^i=0, \,\,\, b^kc^i=0,\,\, k+i>2^s.
\end{equation}

\medskip

The 8th and 9th relations.

\begin{equation*}
\label{eq:Tb3}
(c+x_1+v_s\sum_{i=1}^{s-1}c^{2^s-2^i}x_2^{2^{i-1}})(b+y_1+v_s\sum_{i=1}^{s-1}b^{2^s-2^i}y_2^{2^{i-1}})=v_sTb^{2^s-1}.
\end{equation*}

\begin{proof}
Let $G=G_{38}$. Consider the diagram

\medskip

\begin{equation}
\begin{CD}
\label{eq:DG38b}
B \langle \mathbf{a},\mathbf{c}^2 \rangle @>>> B \langle \mathbf{a},\mathbf{c} \rangle \\
@VV{\rho_{\mu}}V              @VV{\rho_{\beta}}V \\
B \langle \mathbf{a},\mathbf{b},\mathbf{c}^2 \rangle @>>\rho_{\gamma}> BG
\end{CD}
\end{equation}

\medskip

Then the left hand side of our relation is equal to
\begin{align*}
&Tr^*_{\gamma}(u)(b+y_1+v_s\sum_{i=1}^{s-1}b^{2^s-2^i}y_2^{2^{i-1}})               && \text{by transfer formula \eqref{eq:tr}} \\
&=Tr^*_{\gamma}(u) Tr^*_{\beta}(c_1(\lambda'))                                          && \text{by Lemma \ref{lem:LAMBDA38} } \\
&=Tr^*_{\gamma}(u \cdot \rho^*_{\gamma}Tr^*_{\beta}(c_1(\lambda')))                      && \text{by Frobenius reciprocity of the transfer} \\
&=Tr^*_{\gamma}(u \cdot Tr^*_{\mu}(\rho^*_{\mu}(v)))                                  && \text{by the double coset formula and Lemma \ref{lem:LAMBDA38}} \\
&=Tr^*_{\gamma}(Tr^*_{\mu}(\rho^*_{\mu}(uv)))                                      && \text{by Frobenius reciprocity} \\
&=Tr^*_{\gamma}(uv \cdot Tr^*_{\mu}(1))=Tr^*_{\gamma}(uv \cdot v_sc_1^{2^{s}-1}(\mu))      && \text{by the formula for $Tr^*(1)$} \\
&=v_sTb^{2^{s}-1}.                                                                  && \text{by the definitions of $\beta$, $\mu$, and $T$} \\
\end{align*}
\end{proof}

Similarly

\begin{equation*}
\label{eq:Ta3}
(c+y_1+v_s\sum_{i=1}^{s-1}c^{2^s-2^i}y_2^{2^{i-1}})(a+x_1+v_s\sum_{i=1}^{s-1}a^{2^s-2^i}x_2^{2^{i-1}})+v_sTa^{2^s-1}.
\end{equation*}

\begin{proof} Consider the diagram

\medskip

\begin{equation}
\label{eq:DG38a}
\begin{CD}
B\langle \mathbf{a}^2,\mathbf{b},\mathbf{c}^2\rangle           @>>>         B\langle \mathbf{a}^2,\mathbf{b},\mathbf{c}\rangle\\
@VV{\rho_{\lambda^2}}V              @VV{\rho_{\alpha}}V \\
B\langle \mathbf{a},\mathbf{b},\mathbf{c}^2 \rangle             @>>\rho_{\gamma}>         BG
\end{CD}
\end{equation}

\medskip

With this notation the left hand side of the above relation is equal to
\begin{align*}
&Tr^*_{\gamma}(v)  Tr^*_{\alpha}(c_1(\nu'))                                          && \text{by Lemma \ref{lem:NU38} and formula \eqref{eq:tr}}\\
&=Tr^*_{\gamma}(v \cdot \rho^*_{\gamma}Tr^*_{\alpha}(c_1(\nu')))                    && \text{by Frobenius reciprocity of the transfer}\\
&=Tr^*_{\gamma}(v \cdot Tr^*_{\lambda^2}(\rho^*_{\lambda^2}(u)))                        && \text{by the double coset formula and Lemma \ref{lem:NU38}}\\
&=Tr^*_{\gamma}(Tr^*_{\lambda^2}(\rho^*_{\lambda^2}(uv)))                                  && \text{by Frobenius reciprocity}\\
&= Tr^*_{\gamma}(uv \cdot Tr^*_{\lambda^2}(1))=Tr^*_{\gamma}(uv \cdot v_sc_1^{2^{s}-1}(\lambda^2))  && \text{by the formula for $Tr^*(1)$}\\
&=v_sTa^{2^{s}-1}                                                                          && \text{by the definitions of $\alpha$, $\lambda$, and $T$.}\\
\end{align*}
\end{proof}

For the other cases of $G$ the proof is analogous. We just have to apply Lemma \ref{lem:NU39} and Lemma \ref{lem:LAMBDA39} for $G=G_{39}$ or $G_{40},$ and Lemma \ref{lem:LAMBDA41} and Lemma \ref{lem:NU41} for $G_{41}.$

\medskip

The same applies to 11th and 12th relations and may be proved for all four groups simultaneously.
In each case we will arrive at

$$T(a+x_1+v_s\sum_{i=1}^{s-1}a^{2^s-2^i}x_2^{2^{i-1}})=v_sa^{2^s-1}Tr^*(u^2v)$$
or
$$T(b+y_1+v_s\sum_{i=1}^{s-1}b^{2^s-2^i}y_2^{2^{i-1}})=v_sb^{2^s-1}Tr^*(uv^2).$$
Therefore we will need that for the involution $t \in C_2=G/H$ one has by Frobenius reciprocity

\medskip

\noindent i) $ Tr^*(u^2v)=Tr^*(uv(u+tu)-vutu))=Tr^*(uv)x_1-Tr^*(v)x_2$,

\noindent ii) $ Tr^*(uv^2)=Tr^*(uv(v+tv)-uvtv))=Tr^*(uv)y_1-Tr^*(u)y_2.$

\bigskip

Here is details for $G=G_{38}$. Apply again the diagram \eqref{eq:DG38a}.

\begin{align*}
&T(a+x_1+v_s\sum_{i=1}^{s-1}a^{2^s-2^i}x_2^{2^{i-1}})                               && \\
&=Tr^*_{\gamma}(uv) Tr^*_{\alpha}(c_1(\nu'))                                        && \text{by Lemma \ref{lem:NU38} and formula \eqref{eq:tr}}\\
&=Tr^*_{\gamma}(uv \cdot \rho^*_{\gamma}Tr^*_{\alpha}(c_1(\nu')))                    && \text{by Frobenius reciprocity of the transfer}\\
&=Tr^*_{\gamma}(uv \cdot Tr^*_{\lambda^2}(\rho^*_{\lambda^2}(u)))                        && \text{by the double coset formula and Lemma \ref{lem:NU38}}\\
&=Tr^*_{\gamma}(Tr^*_{\lambda^2}(\rho^*_{\lambda^2}(u^2v)))                                             && \text{by Frobenius reciprocity}\\
&= Tr^*_{\gamma}(u^2v \cdot Tr^*_{\lambda^2}(1))=Tr^*_{\gamma}(u^2v \cdot v_sc_1^{2^{s}-1}(\lambda^2))  && \text{by the formula for $Tr^*(1)$}\\
&=v_sT(u^2v)a^{2^{s}-1}                                                                                 && \text{by the definitions of $\alpha$ and $\lambda$}\\
\end{align*}

and the above equality i) gives
$$
T(a+x_1+v_s\sum_{i=1}^{s-1}a^{2^s-2^i}x_2^{2^{i-1}})+v_sa^{2^s-1}Tx_1+v_sa^{2^s-1}x_2(c+y_1+v_s\sum_{i=1}^{s-1}c^{2^s-2^i}y_2^{2^{i-1}})=0.
$$
Then the second summand is zero by the 6th relation. The third summand is equal to $v_sa^{2^s-1}x_2(c+y_1)$ by \eqref{a^kc^l}. This gives the 11th relation.

\medskip

Similarly, applying the diagram \eqref{eq:DG38b} and the above equality ii) we have

$$
T(b+y_1+v_s\sum_{i=1}^{s-1}b^{2^s-2^i}x_2^{2^{i-1}})+b^{2^s-1}Ty_1+b^{2^s-1}y_2(c+x_1+v_s\sum_{i=1}^{s-1}c^{2^s-2^i}x_2^{2^{i-1}})=0.
$$

Again the second summand is zero by the 7th relation and the third summand is equal to $b^{2^s-1}y_2(c+x_1)$ and we get the 12th relation.

\qed

\medskip

The relation $cT=0$ is easy.

\begin{proof}
$cT\equiv cTr^*_{\gamma}(uv)=Tr^*_{\gamma}(uv \gamma^*(c))=Tr^*_{\gamma}(uv \cdot 0)=0$.
\end{proof}

\medskip

Now let us prove the 10th relation.

\begin{proof}
Let $u'=tu$ and $v'=tv$ be as above and $Tr^*=Tr^*_{\gamma}$. Then

\begin{align*}
&Tr^*(uv)+Tr^*(uv')=Tr^*(u(v+v'))=Tr^*(u)Tr^*(v),&&
\end{align*}

\begin{align*}
Tr^*(uv)Tr^*(uv')&=Tr^*(uv(uv'+u'v))  &&\\
&=Tr^*(u^2 vv')+Tr^*(v^2 uu')=Tr^*(u^2)y_2+Tr^*(v^2)x_2.&&
\end{align*}

Also

\medskip

$Tr^*(u^2)=Tr^*(u(u+u')-uu')=Tr^*(u)x_1-Tr^*(1)x_2,$

$Tr^*(v^2)=Tr^*(v(v+v')-vv')=Tr^*(v)y_1-Tr^*(1)y_2.$

\medskip

Now we apply these formulas and take into account $Tr^*(1)x_1=Tr^*(1)y_1=0$. This gives the quadratic equation in $T=Tr^*_{\gamma}(uv)$

\begin{align*}
T^2&=T(c+x_1+v_s\sum_{i=1}^{s-1}c^{2^s-2^i}x_2^{2^{i-1}})(c+y_1+v_s\sum_{i=1}^{s-1}c^{2^s-2^i}y_2^{2^{i-1}})&&\\
&+x_2y_1(c+y_1+v_s\sum_{i=1}^{s-1}c^{2^s-2^i}y_2^{2^{i-1}})+x_1y_2(c+x_1+v_s\sum_{i=1}^{s-1}c^{2^s-2^i}x_2^{2^{i-1}}).&&
\end{align*}
\medskip

Now to get the 10th relation apply $cT=0$.
\end{proof}

\medskip

\medskip

The decompositions for $x_1$ and $y_1$ are the consequences of the formula \eqref{eq:FGL} applied to the determinants of $\nu_!$ and $\lambda_!$.

We need $(x_1x_2)^{2^{2s-2}}=0$ and $(y_1y_2)^{2^{2s-2}}=0$. It follows from the relations ii) of Theorem \ref{thm:G} that moreover we have $(x_1x_2)^{2^s}=(y_1y_2)^{2^s}=0$: decompositions $x_1^{2^s}=(ac)^{2^{s-1}}$ and $y_1^{2^s}=(bc)^{2^{s-1}}$ imply

$$x_1^{2^s}a=x_1^{2^s}c=x_1^{2^s}b^2=y_1^{2^s}b=y_1^{2^s}c=y_1^{2^s}a^2=0$$

since $a^2c=ac^2$, $b^2c=bc^2$ and $a^{2^s}=b^{2^s}=c^{2^s}=0.$

\medskip

That is, all the terms of the above decomposition of $x_2^{2^s}$ annihilate $x_1^{2^s}$. Similarly for $y_1$ and $y_2$.
It is also clear that for computing the Euler classes of the determinants (see Section 3) in each case we need only the initial fragment of the formal group law $F(x,y)=x+y+v_s(xy)^{2^{s-1}}.$

For instance consider $G_{38}$: $det(\lambda_!)=\alpha \beta \gamma$, $C(\lambda_!)=1+y_1+y_2$, $a=e(\alpha)$, $b=e(\beta)$, $c=e(\gamma)$,  hence $e(det\lambda_!)=F(a,F(b,c))=a+b+c+v_s(ab+ac+bc)^{2^{s-1}}.$

For $G_{41}$ we have different $a$ and $b$, but  $det(\lambda_!)=\beta \gamma=\alpha(\alpha \beta)\gamma$, hence $e(det(\lambda_!))=F(a,F(b,c))$ as for $G_{38}$.

\qed

%
%
%
%

\medskip

\section{Invariants}

Let $\alpha$, $\beta$, $\lambda$, $\nu$, $\lambda_!$, $\nu_!$ be as above. We need the action of the involution $t\in G/H=C_2$ on $K(s)^*(BH)$. For simplicity we will ignore the powers of $v_s$. Also we will denote the restrictions of the generators of Theorem \ref{thm:G} to $K(s)^*(BH)$ with the same symbols but with bars.

\begin{lemma}
\label{lem:action}
Let $G=G_{39},G_{40}$. Then

$t(u)=u+u^{2^s}+v^{2^s}+(uv)^{2^{s-1}}+u^{2^{s-1}(1+2^s)}+u^{2^{s-1}}v^{2^{2s-1}};$

$t(v)=v+v^{2^s}+v^{2^{s-1}(1+2^s)}.$

\end{lemma}

\begin{proof}
We need the action of the involution on bundles in Section 3.
$$t(\lambda)=\lambda^3\nu^2=\lambda(\lambda \nu)^2,\,\,\, t(\nu)=\nu^3=\nu(\nu)^2.$$
Note that the initial segment of the formal group law suffices. Namely as $\lambda^4=\nu^4=1$ we can apply the formula $F(y,z)=y+z+(yz)^{2^{s-1}}$ modulo $z^{2^{2(s-1)}}$ (see \cite{BP1}, Lemma 5.3)

\medskip

\begin{align*}
 & t(u)=F(u,F(u^{2^s},v^{2^s}))&&\text{}\\
 = &u+(u^{2^s}+v^{2^s}+(uv)^{2^{2s-1}})+u^{2^{s-1}}(u^{2^s}+v^{2^s})^{2^{s-1}}&&\text{}\\
 = &u+u^{2^s}+v^{2^s}+(uv)^{2^{2s-1}}+u^{2^{s-1}(1+2^s)}+u^{2^{s-1}}v^{2^{2s-1}}.&&\
\end{align*}

\medskip

Similarly for $v$
$$t(v)=F(v,v^{2^s})=v+v^{2^s}+v^{2^{s-1}(1+2^s)}.$$

\end{proof}

We recall that $K(s)^*(BH)=\mathbb{F}_p[v_s,v_s^{-1}][u,v]/(u^{4^s},v^{4^s}).$ This is because of the K\"{u}nneth isomorphisms,
$K(s)^*(BC_4 \times BC_4) = K(s)^*(BC_4) \otimes K(s)^*(BC_4)$,
whereas Morava $K$-theories for cyclic groups are the truncated polynomials \cite{RAV}.
So in particular $K(s)^*(BC_4)=K(s)^*[u]/u^{4^s}$.

Then we have by definition, where bar is defined as in the first paragraph of this section, that

\medskip

$\bar{a}=F(u^{2^s},v^{2^s})=u^{2^s}+v^{2^s}+(uv)^{2^{2s-1}};$

$\bar{b}=v^{2^s};$

$\bar{x_1}=u+t(u)=u^{2^s}+v^{2^s}+(uv)^{2^{2s-1}}+u^{2^{s-1}(1+2^s)}+u^{2^{s-1}}v^{2^{2s-1}};$

$\bar{x_2}=ut(u)=u(u+u^{2^s}+v^{2^s}+(uv)^{2^{2s-1}}+u^{2^{s-1}(1+2^s)}+u^{2^{s-1}}v^{2^{2s-1}});$

$\bar{y_1}=v+t(v)=v^{2^s}+v^{2^{s-1}(1+2^s)};$

$\bar{y_2}=vt(v)=v(v+v^{2^s}+v^{2^{s-1}(1+2^s)});$

$\bar{T}=uv+t(uv)$

$=uv+(u+u^{2^s}+v^{2^s}+(uv)^{2^{2s-1}}+u^{2^{s-1}(1+2^s)}+u^{2^{s-1}}v^{2^{2s-1}})(v+v^{2^s}+v^{2^{s-1}(1+2^s)}).$

\medskip

Note that as $u^{2^{2s}}=v^{2^{2s}}=0$ one has
$$
\bar{x}_1^{2^s}=\bar{y}_1^{2^s}=\bar{x}_2^{2^{2s-1}}=\bar{y}_2^{2^{2s-1}}=0.
$$

\medskip

To describe all invariants we need the following

\begin{lemma}
\label{lem:free4module}
i) Let $G=G_{39},G_{40}.$ Then $K(s)^*(BH)$ is a free $K(s)^*(\bar{x_2},\bar{y_2})/(\bar{x_2}^{2^{2s-1}},\bar{y_2}^{2^{2s-1}})$ module generated by $1, u, v, uv;$

\medskip

ii) $K(s)^*$-rank of $K(s)^*(BH)^{\mathcal{C}_2}$ is $16^s/2+4^s/2$ and a basis is

\medskip

1) $\bar{x_2}^i \bar{y_2}^j uv$,
2) $\bar{x_2}^i \bar{y_2}^j u$,\,\,
$i,j \geq 2^{2s-1}-2^{s-1}$ in 1), 2);

3) $\bar{x_2}^i \bar{y_2}^j v$, $j\geq 2^{2s-1}-2^{s-1} $;

4) $\bar{x_2}^i \bar{y_2}^j u+
\bar{x_2}^i\bar{y_2}^{j-2^{s-1}}(\bar{x_2}^{2^{s-1}}+\bar{y_2}^{2^{s-1}}\sum_{k=1}^{s}\bar{y_2}^{2^{2s-1}-2^{2s-k}+2^{s-k}-2^{s-1}}+
\bar{x_2}^{2^{2s-2}}\bar{y_2}^{2^{2s-2}})v$,

$i<2^{2s-1}-2^{s-1}$, $j\geq 2^{s-1}$;

\medskip

5) $\bar{x_2}^i \bar{y_2}^ju
+\bar{x_2}^i\bar{y_2}^j(\sum_{k=1}^{s}\bar{y_2}^{2^{2s-1}-2^{2s-k}+2^{s-k}-2^{s-1}})v+\bar{x_2}^{i-2^{s-1}}\bar{y_2}^{j+2^{2s-1}-2^{s-1}}uv$,

$i\geq 2^{2s-1}-2^{s-1}$, $j<2^{s-1}$;

\medskip

6) $\bar{x_2}^i\bar{y_2}^ju+\bar{x_2}^i\bar{y_2}^j(\sum_{k=1}^{s}\bar{y_2}^{2^{2s-1}-2^{2s-k}+2^{s-k}-2^{s-1}})v,$

$i\geq 2^{2s-1}-2^{s-1}$, $2^{s-1}\leq j<2^{2s-1}-2^{s-1}$.

\medskip

7) $\bar{x_2}^i\bar{y_2}^j$;

\medskip

$i,j<2^{2s-1}$ in 1)-7).

\medskip

\medskip

iii) There is a decomposition of $K(s)^*(BH)$ into free and trivial $\mathcal{C}_2$ modules, such that a basis for the trivial module is
$\bar{x_2}^i\bar{y_2}^j,$ $\bar{a}\bar{x_2}^i\bar{y_2}^j$ $i<2^{s}, j<2^{s-1}$.

\medskip

iv) $K(s)^*(BG)$ is generated by  $c,a,b,x_2,y_2,T$ as a $K(s)^*$-algebra.

\end{lemma}

\medskip

\begin{proof}

i) It suffices to prove that $\{\bar{x_2}^i\bar{y_2}^i\}\times\{1,u,v,uv\}$, $i,j<2^{2s-1}$
is a $K(s)^*$ generating set. Clearly it will be a $K(s)^*$ basis by its number of elements. Any polynomial in $u,v$ can be uniquely written as $f_0+f_1u+f_2v+f_3uv$, $f_i=f_i(\bar{x_2},\bar{y_2})$ as follows. Because $u^2=u\bar{x_1}-\bar{x_2}$ and $v^2=v\bar{y_1}-\bar{y_2}$ any polynomial in $u,v$ can be uniquely written as  $g_0+g_1u+g_2v+g_3uv$ where $g_i=g_i(\bar{x_1},\bar{y_1},\bar{x_2},\bar{y_2}).$ But

\begin{equation}
\label{eq:decy1}
\bar{y_1}=\bar{y_2}^{2^{s-1}}.
\end{equation}

This follows from the decomposition of $\bar{y_1}$ of Theorem \ref{thm:G} as $\bar{y_2}^{2^{2s-1}}=0:$

$\bar{y_1}=\bar{y_2}^{2^{s-1}}+\bar{y_1}^{2^{s-1}}\bar{y_2}^{2^{2s-2}}=\bar{y_2}^{2^{s-1}}\, modulo \,\,\, \bar{y_2}^{2^{2s-2}}\bar{y_2}^{2^{2s-2}}=\bar{y_2}^{2^{2s-1}}.$

Similarly the decomposition of $\bar{x}_1$ of Theorem \ref{thm:G} implies

\begin{align*}
&\bar{x_1}=\bar{b}+\bar{x_2}^{2^{s-1}}+\bar{x_1}^{2^{s-1}}\bar{x_2}^{2^{2s-2}} && \\
 =&\bar{b}+\bar{x_2}^{2^{s-1}}+\bar{b}^{2^{s-1}}\bar{x_2}^{2^{2s-2}} \,modulo \,\,\,\bar{x_2}^{2^{2s-2}}\bar{x_2}^{2^{2s-2}}  && \\
 =&\bar{b}+\bar{x_2}^{2^{s-1}}+\bar{b}^{2^{s-1}}\bar{x_2}^{2^{2s-2}}&& \text{(as $\bar{x_2}^{2^{2s-1}}=0$).}\
\end{align*}

 Then
 \begin{equation}
\label{eq:u^2^s}
u^{2^s}=u\bar{x_1}^{2^s-1}+\sum_{i=1}^s\bar{x_1}^{2^s-2^i}\bar{x_2}^{2^{i-1}}
\end{equation}
follows inductively from $u^2=u\bar{x}_1+\bar{x}_2.$ We can replace $u,\bar{x}_i$ by $v,\bar{y}_i$ in \eqref{eq:u^2^s} and get (by inductive argument again) from $v^2=v\bar{y}_1+\bar{y}_2$ and \eqref{eq:decy1}

\medskip

\begin{equation}
\label{eq:decb}
\bar{b}=v^{2^s}=\sum_{i=1}^{s}\bar{y_2}^{2^{2s-1}+2^{s-i}-2^{2s-i}}+v\bar{y_2}^{2^{2s-1}-2^{s-1}}.
\end{equation}


Then by Theorem \ref{thm:G} $\bar{b}^2=\bar{y_2}^{2^s}$ and we get for $\bar{x_1}$

\begin{equation}
\label{eq:decx1}
\bar{x_1}=\bar{x_2}^{2^{s-1}}+\bar{y_2}^{2^{s-1}}\sum_{i=1}^{s}\bar{y_2}^{2^{2s-1}-2^{2s-i}+2^{s-i}-2^{s-1}}+
\bar{x_2}^{2^{2s-2}}\bar{y_2}^{2^{2s-2}}+v\bar{y_2}^{2^{2s-1}-2^{s-1}}.
\end{equation}

ii) First let us write $\bar{T}$ in our basis. Note $t(u)t(v)=(\bar{x_1}-u)(\bar{y_1}-v)=\bar{x_1}\bar{y_1}-\bar{x_1}v-\bar{y_1}u+uv$ and we have $\bar{T}=uv+t(uv)=\bar{x_1}\bar{y_1}-\bar{x_1}v-\bar{y_1}u$. But

$$
\bar{x_1}v=(\bar{x_2}^{2^{s-1}}+\bar{y_2}^{2^{s-1}}\sum_{i=1}^{s}\bar{y_2}^{2^{2s-1}-2^{2s-i}+2^{s-i}-2^{s-1}}+
\bar{x_2}^{2^{2s-2}}\bar{y_2}^{2^{2s-2}})v+\bar{y_2}^{2^{2s-1}-2^{s-1}+1}
$$
as $v^2=v\bar{y_2}^{2^{s-1}}+\bar{y_2}$ and $\bar{y_2}^{2^{2s-1}}=0$. Then \eqref{eq:decx1} and \eqref{eq:decy1} imply

\begin{equation}
\label{eq:decT}
\begin{split}
&\bar{T}=(\bar{x_2}^{2^{s-1}}+\bar{y_2}^{2^{s-1}}\sum_{i=1}^{s}\bar{y_2}^{2^{2s-1}-2^{2s-i}+2^{s-i}-2^{s-1}}+\bar{x_2}^{2^{2s-2}}\bar{y_2}^{2^{2s-2}})\bar{y_2}^{2^{s-1}}+
\bar{y_2}^{2^{2s-1}-2^{s-1}+1}\\
+&(\bar{x_2}^{2^{s-1}}+\bar{y_2}^{2^{s-1}}\sum_{i=1}^{s}\bar{y_2}^{2^{2s-1}-2^{2s-i}+2^{s-i}-2^{s-1}}+\bar{x_2}^{2^{2s-2}}\bar{y_2}^{2^{2s-2}})v \\
+&\bar{y_2}^{2^{s-1}}u.
\end{split}
\end{equation}

\medskip

\medskip

Now let $g=f_0+f_1u+f_2v+f_3uv$, $f_i=f_i(\bar{x_2},\bar{y_2})$ be an invariant, that is, $g \in Ker(1+t).$ Then

\begin{equation*}
\label{eq:f1f2f3}
f_1\bar{x_1}+f_2\bar{y_1}+f_3\bar{T}=0.
\end{equation*}

Taking into account the decompositions \eqref{eq:decx1}, \eqref{eq:decy1} and \eqref{eq:decT} we have

\begin{equation}
\label{eq:system}
\begin{split}
&f_3\bar{y_2}^{2^{s-1}}=0, \,\,\text{therefore}\\
&f_3\bar{x_2}^{2^{s-1}}=f_1\bar{y_2}^{2^{2s-1}-2^{s-1}},\\
&f_2\bar{y_2}^{2^{s-1}}=f_1( \bar{x_2}^{2^{s-1}}+\bar{y_2}^{2^{s-1}}\sum_{i=1}^{s}\bar{y_2}^{2^{2s-1}-2^{2s-i}+2^{s-i}-2^{s-1}}+
\bar{x_2}^{2^{2s-2}}\bar{y_2}^{2^{2s-2}}).\\
\end{split}
\end{equation}

\medskip

Now by the third equation of \eqref{eq:system} we have two possible cases
\newline a) $f_1$ has a factor $\bar{y_2}^{2^{s-1}}$ and we can restore $f_2$ modulo summands of type 3). Also $f_3=0$ modulo summands of type 1). Therefore  $g$ is decomposable into sum of elements of types 1), 3), 4), 7);
\newline b) $f_1$ annihilates $\bar{x_2}^{2^{s-1}}$, that is $f_1$ has a factor $\bar{x_2}^{2^{2s-1}-2^{s-1}}$. Hence right hand side of the third equation of \eqref{eq:system} has a factor $\bar{y_2}^{2^{s-1}}$ and we can restore $f_2$ modulo summands of type 3). Then by the second equation of \eqref{eq:system} we can restore $f_3$ modulo summands of type 1). Hence $g$ is decomposable into elements of types 1), 3), 5)(or 6)) and 7).
If $f_1$ annihilates $\bar{x_2}^{2^{s-1}}$ and $\bar{y_2}^{2^{s-1}}$ then by \eqref{eq:system} $g$ is decomposable into the monomials 1), 2), 3) and 7) of Lemma \ref{lem:free4module}.

 By Lemma \ref{lem:free4module} i) the elements 1)-7) are independent, therefore they form a basis ($16^s/2+4^s/2$ elements in total) for the invariants.

\medskip

iii) Consider now the decomposition $[K(s)^*(H)]^{\mathcal{C}_2}=(\mathcal{F})^{\mathcal{C}_2}+\mathcal{T}$, corresponding to the decomposition of  $K(s)^*(H)$ into free and trivial $\mathcal{C}_2$-modules.
\newline Clearly the composition $\rho^*Tr^*=1+t$ is onto $(\mathcal{F})^{\mathcal{C}_2}$. Also i) and the Fr\"{o}benius reciprocity of the transfer implies $Im Tr^* \subset A$, where $A$ is the subalgebra in $K(s)^*(BG)$ generated by $c,a,b,x_2,y_2,T$. Here we use the formulas \eqref{Tru}, \eqref{Trv} for $Tr^*(u)$, $Tr^*(v)$ respectively and the decompositions of $x_1$ and $y_1$ of Theorem \ref{thm:G}. We have to check whether the invariants in $\mathcal{T}$ are also covered by $\rho^*$. Let $m=\chi_s(\mathcal{F}^{\mathcal{C}_2})$ and $n=\chi_s(\mathcal{T})$ be the $K(s)^*$-Euler characteristics. Then
$$\chi_s(H)=2m+n, \,\,\,\chi_s(H)^{\mathcal{C}_2}=m+n.$$
Clearly $2m+n=16^s$ and by ii) $m+n=16^s/2+4^s/2$. It follows that $m=16^s/2-4^s/2$ and $n=4^s.$

Let us consider the invariants modulo $Im(1+t)$ and denote by $S_k$, $k=1,\cdots,7$, the set of invariants of type k) of Lemma \ref{lem:free4module} ii) modulo $Im(1+t).$
First note that

\begin{equation}
\label{eq:xy}
S_7 \subset S_8, \text{where} \,\,S_8=\{\bar{x_2}^i\bar{y_2}^j, i<2^s, j<2^{s-1}\},
\end{equation}
\newline that is except $\{\bar{x_2}^i\bar{y_2}^j, i<2^s, j<2^{s-1}\}$ all monomials of Lemma \ref{lem:free4module} of type 7) $\equiv 0.$
This is because by definition $\bar{x}_1,\bar{y}_1 \equiv 0;$ By \eqref{eq:decy1} $\bar{y_2}^{2^{s-1}}\equiv\bar{y}_1$ ; The decomposition of $\bar{x}_1$ after \eqref{eq:decy1} implies $\bar{x_2}^{2^{s}}\equiv \bar{b}^2\equiv 0$ as by Theorem \ref{thm:G} $\bar{b}^2\equiv \bar{y_2}^{2^s}\equiv \bar{y_1}^2\equiv 0$.

Then by \eqref{eq:decx1} one has

\begin{equation}
\label{eq:vy}
\bar{y_2}^{2^{2s-1}-2^{s-1}}v \equiv \bar{x_2}^{2^{s-1}},\text{\,\,\,hence\,\,\,}S_3 \subset S_8.
\end{equation}
Also it follows all monomials of type 1) $\equiv 0$, that is $S_1=\{0\}$ as $\bar{x}_2^{2^{2s-1}}=0.$

\medskip

\eqref{eq:decT} implies

$$\bar{y_2}^{2^{s-1}}u=\bar{y_2}^{2^{s-1}}\sum_{i=1}^{s}\bar{y_2}^{2^{2s-1}-2^{2s-i}+2^{s-i}-2^{s-1}}v,\,\, modulo\,\, \bar{x_2}^{2^{s-1}},$$
hence $S_2 \subset S_3$. Thus by \eqref{eq:vy} $S_2 \subset S_8.$

\medskip

Recall by definition $\bar{T}\in Im (1+t).$ Also the first line in the decomposition \eqref{eq:decT} $\in Im(1+t)$ as it has factor $\bar{y_2}^{2^{s-1}}=\bar{y}_1$. Therefore if one denotes by $\breve{T}$ the sum of second and third lines in \eqref{eq:decT} one has
$$
0\equiv \breve{T}=(\bar{x_2}^{2^{s-1}}+\bar{y_2}^{2^{s-1}}\sum_{i=1}^{s}\bar{y_2}^{2^{2s-1}-2^{2s-i}+2^{s-i}-2^{s-1}}+
\bar{x_2}^{2^{2s-2}}\bar{y_2}^{2^{2s-2}})v+\bar{y_2}^{2^{s-1}}u.$$
Then the elements of $S_4$ are $\bar{x_2}^i\bar{y_2}^{j-2^{s-1}}\breve{T}\equiv 0;$
\newline
The same argument implies $S_6=\{0\}.$  Now only the elements of $S_5$ remain. It suffices to prove that the element of $S_5$ obtained for $i={2^{2s-1}}-2^{s-1}$ and $j=0$ corresponds to $\bar{a}$, that is

$$\bar{a}\equiv\bar{x_2}^{2^{2s-1}-2^{s-1}}u
+\bar{x_2}^{2^{2s-1}-2^{s-1}}(\sum_{k=1}^{s}\bar{y_2}^{2^{2s-1}-2^{2s-k}+2^{s-k}-2^{s-1}})v+\bar{x_2}^{2^{2s-1}-2^{s}}\bar{y_2}^{2^{2s-1}-2^{s-1}}uv.
$$

First note

\begin{equation}
\label{eq:amod(1+t)}
\bar{a}\equiv u\bar{x_1}^{2^s-1}:
\end{equation}

Recall $\bar{a}=u^{2^s}+v^{2^s}+(uv)^{2^{2s-1}}$. Then \eqref{eq:u^2^s} implies
$u^{2^s}\equiv u\bar{x_1}^{2^s-1}+\bar{x_2}^{2^{s-1}}$. Recall $\bar{b}=v^{2^s}$ and by \eqref{eq:decb} and \eqref{eq:decx1} $\bar{x_2}^{2^{s-1}}+b\equiv 0.$ It suffices to see $(uv)^{2^{2s-1}}\equiv 0.$ \eqref{eq:u^2^s} also implies $u^{2^{2s-1}}=\bar{x_2}^{2^{2s-2}}$. Similarly $v^{2^{2s-1}}=\bar{y_2}^{2^{2s-2}}$. Then recall by \eqref{eq:decy1} $\bar{y_2}^{2^{s-1}}\equiv 0$ hence  $(uv)^{2^{2s-1}}=\bar{x_2}^{2^{2s-2}}\bar{y_2}^{2^{2s-2}}\equiv 0$.

\medskip

Now it follows (see the proof of \eqref{eq:decx1}) from $\bar{x_1}=\bar{b}+\bar{x_2}^{2^{s-1}}+\bar{b}^{2^{s-1}}\bar{x_2}^{2^{2s-2}}$ and $\bar{b}^2=\bar{y_2}^{2^s}$ that

\begin{align*}
u\bar{x_1}^{2^s-1}=&u\bar{x_1}\bar{x_1}^{2^s-2}=ux_1(\bar{x_2}^{2^s}+\bar{y_2}^{2^s})^{2^{s-1}-1}&&\\
=&(u \bar{b}+u\bar{x_2}^{2^{s-1}})(\bar{x_2}^{2^s}+\bar{y_2}^{2^s})^{2^{s-1}-1}
+u(\bar{x_2}\bar{y_2})^{2^{2s-2}}(\bar{x_2}^{2^s}+\bar{y_2}^{2^s})^{2^{s-1}-1}&&\\
=&(u \bar{b}+u\bar{x_2}^{2^{s-1}})(\bar{x_2}^{2^s}+\bar{y_2}^{2^s})^{2^{s-1}-1}.&&\
\end{align*}

Now let

$$(u \bar{b}+u\bar{x_2}^{2^{s-1}})(\bar{x_2}^{2^s}+\bar{y_2}^{2^s})^{2^{s-1}-1}=f_0+f_1u+f_2v+f_3uv,$$
where $f_i$ are polynomials in $\bar{x_2},\bar{y_2}$. Let us prove that $f_3=\bar{x_2}^{2^{2s-1}-2^{s}}\bar{y_2}^{2^{2s-1}-2^{s-1}}.$
Then we will restore $f_1$ and $f_2$ as above and complete the proof.

To get $f_3$ note that by \eqref{eq:decb} only the summand $v\bar{y_2}^{2^{2s-1}-2^{s-1}}$ in the decomposition of $\bar{b}$ is relevant. This gives
$$uv\bar{y_2}^{2^{2s-1}-2^{s-1}}(\bar{x_2}^{2^s}+\bar{y_2}^{2^s})^{2^{s-1}-1}=uv\bar{x_2}^{2^{2s-1}-2^s}\bar{y_2}^{2^{2s-1}-2^{s-1}}$$
and proves our claim for $f_3.$ Then by \eqref{eq:system}
$$f_1\equiv \bar{x_2}^{2^{2s-1}-2^{s-1}},\,\,\, f_2\equiv \bar{x_2}^{2^{2s-1}-2^{s-1}}\sum_{j=1}^{s}\bar{y_2}^{2^{2s-1}-2^{2s-j}+2^{s-j}-2^{s-1}}$$
modulo $\bar{x_2}^k\bar{y_2}^l$, $l\geq 2^{s-1}$, which are elements of $Im(1+t)$. Hence $a \equiv f_1u+f_2v+f_3uv\equiv$ the element of $S_5$ with $i=2^{2s-1}-2^{s-1}$ and $j=0$. After multiplying by $\bar{x_2}^i\bar{y_2}^j$, $i<2^s, j<2^{s-1}$ we get  $a\bar{x_2}^i\bar{y_2}^j$ modulo $Im(1+t)$ (still elements of $S_5$).

\medskip

iv) Consider the Serre spectral sequence for the extension $1 \to H \to G \to G/H \to 1$.  Then as $G/H=\mathcal{C}_2$
$$E_2=H^*(\mathcal{C}_2,K(s)^*(BH))\cong \mathcal{F}^{\mathcal{C}_2}\oplus \mathcal{T} \otimes H^*(\mathcal{C}_2,\mathbb{F}_2),$$
where $K(s)^*(BH)=\mathcal{F}\oplus \mathcal{T}$ is the decomposition into free and trivial modules.
As in iii) let $A$ be the subalgebra of $K(s)^*(BG)$ generated by $c,a,b,x_2,y_2,T$. Then iii) says that the restriction
$\rho^*: A \to [K(s)^*(BH)]^{\mathcal{C}_2}$ is onto. $\rho^*(c)=0$ by definition. Hence all invariants are permanent cycles and there is only one differential $d_{2^{s+1}-1}(t)=v_st^{2^{s+1}}$. Since $t^2$ is represented by $c$ one obtains that $c,a,b,x_2,y_2,T$ are $K(s)^*$-algebra generators of $K(s)^*(BG)$ and its $K(s)^*$-rank equals $\chi_{s}(\mathcal{F}^{C_2}\oplus \mathcal{T} \otimes \mathbb{F}_2[c]/c^{2^s})= (16^s-4^s)/2+4^s2^s=16^s/2+8^s-4^s/2$ as it is already known from \cite{SCH1}.

This completes the proof.

\end{proof}

\bigskip

\begin{lemma}
\label{lem:actiong38}
Let $G=G_{38}, G_{41}$. Then

$t(u)=u+v^{2^s}+u^{2^{s-1}}v^{2^{2s-1}};$

$t(v)=v+\bar{b}+v^{2^{s-1}}\bar{b}^{2^{s-1}}.$
\end{lemma}

\begin{proof}

Recall the action of the involution on bundles in Section 3. For $G_{38}$
$$t(\lambda)=\lambda\mu,\,\,\, t(\nu)=\lambda^2\nu$$
and for $G_{41}$
$$
t(\lambda)=\lambda \rho^*(\alpha\beta),\,\,\, t(\nu)=\lambda^2\nu.
$$
Since $\lambda^4=\mu^2=(\alpha\beta)^2=1$ again we need only the initial segment of the formal group law.
\end{proof}

Recall that $H$ is isomorphic to $C_4\times C_2\times C_2$ for $G=G_{38}$ and to $C_4\times C_4$ for $G=G_{41}$.
Again by the K\"{u}nneth isomorphism we have that as a $K(s)^*$-algebra $K(s)^*(BG_{41})=K(s)^*[u,v]/(u^{4^s},v^{4^s})$
and $K(s)^*(BG_{38})=K(s)^*[u,v,w]/(u^{4^s},v^{2^s},w^{2^s})$, where $w=c_1(\mu)$ is invariant under action of $G/H=\mathcal{C}_2$.

Then we have by definition

\medskip

$\bar{x_1}=u+t(u)=v^{2^s}+u^{2^{s-1}}v^{2^{2s-1}};$

$\bar{x_2}=ut(u)=u(u+v^{2^s}+u^{2^{s-1}}v^{2^{2s-1}});$

$\bar{y_1}=v+t(v)=\bar{b}+v^{2^{s-1}}\bar{b}^{2^{s-1}};$

$\bar{y_2}=vt(v)=v(v+\bar{b}+v^{2^{s-1}}\bar{b}^{2^{s-1}});$

and as $\bar{T}=uv+t(uv)=uv+(\bar{x_1}+u)(\bar{y_1}+v)$ we have

\begin{equation}
\label{eq:38T}
\bar{T}=\bar{x_1}\bar{y_1}+\bar{x_1}v+\bar{y_1}u.
\end{equation}

\medskip

Now to describe all invariants and see that $K(s)^*(BG)$ restricts onto $K(s)^*(BH)^{\mathcal{C}_2}$, we turn to the following

\begin{lemma}
\label{lem:H38free4module}
i) Let $G=G_{38},G_{41}$ and let $x^{\omega}=\bar{x_1}^i\bar{y_1}^j\bar{x_2}^k\bar{y_2}^l$, $i,j<2^s,\,\,k,l<2^{s-1}$. Then the set $x^{\omega},x^{\omega}u,x^{\omega}v,x^{\omega}uv$ is a $K(s)^*$ basis in $K(s)^*(BH)$.

\medskip

ii) $K(s)^*$-rank of $K(s)^*(BH)^{\mathcal{C}_2}$ is $16^s/2+4^s/2$ and a basis is

1) $\bar{x_1}^i\bar{y_1}^j\bar{x_2}^k\bar{y_2}^l$, $i,j<2^s,\,k,l<2^{s-1}$;

2) $\bar{x_1}^{2^s-1}\bar{y_1}^{2^s-1}\bar{x_2}^k\bar{y_2}^luv$, $k,l<2^{s-1}$;

3) $\bar{x_1}^{2^s-1}\bar{y_1}^i\bar{x_2}^k\bar{y_2}^lu$,\,\,\, $\bar{x_1}^i\bar{y_1}^{2^s-1}\bar{x_2}^k\bar{y_2}^lv$,\,\,\, $i<2^s$, $k,l<2^{s-1}$;

4) $\bar{x_1}^i\bar{y_1}^j\bar{x_2}^k\bar{y_2}^l u + \bar{x_1}^{i+1}\bar{y_1}^{j-1}\bar{x_2}^k\bar{y_2}^l v$,\,\, $i,j-1<2^s-1,\,k,l<2^{s-1}$.

iii) The set $x_2^iy_2^j$, $ax_2^iy_2^j$, $bx_2^iy_2^j$, $abx_2^iy_2^j$,\,\,$i,j<2^{s-1}$ restricts to a $K(s)^*$ basis in the trivial summand $\mathcal{T}$ of the $\mathcal{C}_2$ module $K(s)^*(H)$.

iv) $K(s)^*(BG)$ is generated by  $c,a,b,x_2,y_2,T$ as a $K(s)^*$-algebra.

\end{lemma}

\begin{proof}i) As in Lemma \ref{lem:free4module}.i) any polynomial $g(u,v)$ can be uniquely written as  $g_0+g_1u+g_2v+g_3uv$ where
 $g_i=g_i(\bar{x_1},\bar{y_1},\bar{x_2},\bar{y_2}).$ Then by Theorem \ref{thm:G}  $\bar{x_2}^{2^{s-1}}$, $\bar{y_2}^{2^{s-1}}$ can be expressed by $\bar{x_1},\bar{y_1},\bar{a},\bar{b}$. Now
 $$\bar{x_1}=u+F(u,\bar{a})=u+u+\bar{a}+u^{2^{s-1}}\bar{a}^{2^{s-1}}=\bar{a}+u^{2^{s-1}}\bar{x_1}^{2^{s-1}}$$
and we get
\begin{equation}
\label{eq:38deca}
\bar{a}=\bar{x_1}+\bar{x_1}^{2^s-1}u+\sum_{i=1}^{s-1}\bar{x_1}^{2^{s}-2^i}\bar{x_2}^{2^{i-1}}.
\end{equation}
Similarly $\bar{y_1}=\bar{b}+v^{2^{s-1}}\bar{y_1}^{2^{s-1}}$ implies
\begin{equation}
\label{eq:38decb}
\bar{b}=\bar{y_1}+\bar{y_1}^{2^s-1}v+\sum_{i=1}^{s-1}\bar{y_1}^{2^{s}-2^i}\bar{y_2}^{2^{i-1}}.
\end{equation}
Because of nilpotence of $\bar{x_1}$ and $\bar{y_1}$ substituting $\bar{a}$ and $\bar{b}$ in $g(u,v)$ we arrive at i) after finite number of steps.

\medskip

ii) Let
$$g=f_0+f_1u+f_2v+f_3uv,\,\,\,f_i=f_i(\bar{x_1},\bar{y_1}\bar{x_2},\bar{y_2})$$
be an invariant, that is, $g \in Ker(1+t).$ Then by \eqref{eq:38T}
$$f_1\bar{x_1}+f_2\bar{y_1}+f_3\bar{T}=f_1\bar{x_1}+f_2\bar{y_1}+f_3(\bar{x_1}\bar{y_1}+\bar{x_1}v+\bar{y_1}u)=0$$
and we get

\begin{equation}
\label{eq:38f1f2f3}
f_3\bar{x_1}=f_3\bar{y_1}=0;\,\,f_1\bar{x_1}=f_2\bar{y_1}.
\end{equation}

Now i) implies  ii).

\bigskip

iii) Let us look at invariants modulo $1+t$: The invariants of Lemma \ref{lem:H38free4module} 1) with $i,j>0$ are zero as $\bar{x_1},\bar{y_1}\in Im(1+t);$
Invariants of 4) are all zero as by \eqref{eq:38T}

\begin{equation}
\label{eq:38xv-yu}
\bar{x_1}v+\bar{y_1}u\in Im(1+t).
\end{equation}
By the same argument invariants of 3) are zero except $\bar{x_1}$ or $\bar{y_1}$ are omitted, that is except $\bar{x_1}^{2^s-1}\bar{x_2}^k\bar{y_2}^lu$,\,$\bar{y_1}^{2^s-1}\bar{x_2}^k\bar{y_2}^lv$.
This completes the basis invariants of $\mathcal{F}^{\mathcal{C}_2}$ because of the total number $16^s/2-4^s/2$ (see the proof of Lemma \ref{lem:free4module}. iii).

Thus the invariants of 2) are all nonzero. Then we have that the basis invariants corresponding to $\mathcal{T}$ are as follows:
\newline$\bar{x_2}^i\bar{y_2}^j$ is resticted $x_2^iy_2^j$, \,\,$i,j<2^{s-1}$;
and modulo $1+t$
\newline$\bar{x_1}^{2^s-1}\bar{x_2}^i\bar{y_2}^ju$ is restricted $ax_2^iy_2^j$ by \eqref{eq:38deca};
\newline$\bar{y_1}^{2^s-1}\bar{x_2}^i\bar{y_2}^jv$ is restricted $bx_2^iy_2^j$ by \eqref{eq:38decb};
\newline $\bar{x_1}^{2^s-1}\bar{y_1}^{2^s-1}\bar{x_2}^i\bar{y_2}^juv$ is restricted $abx_2^iy_2^j$ by \eqref{eq:38deca}, \eqref{eq:38decb} and \eqref{eq:38xv-yu}.

iv) This is the consequence of the arguments similar to that of Lemma \ref{lem:free4module} iv).

\end{proof}

\medskip

\begin{remark}
Of course there are alternative bases for $K(s)^*(BH).$ For instance Lemma \ref{lem:free4module}.i) is also true for $G=G_{41};$ For $G=G_{38},G_{41}$ an alternative $K^*(s)$-basis is $X^{\omega},X^{\omega}u,X^{\omega}v,X^{\omega}uv$, where $X^{\omega}=\bar{x_1}^i\bar{x_2}^j\bar{y_2}^k,$ $i<2^s,$ $j<2^{2s-1},$ $k<2^{s-1}.$
\end{remark}

\section*{End of the proof}

By Lemma \ref{lem:free4module} iv) and Lemma \ref{lem:H38free4module} iv) we have that $c,a,b,x_2,y_2,T$ is a complete set of $K(s)^*$ algebra generators of $K(s)^*(BG)$. Now we want to verify that for all our groups the defining relations of Theorem \ref{thm:G} give us a ring of Euler characteristic already computed in \cite{SCH1} $\chi_{s,2}=16^s/2+8^s-4^s/2$.
For each of our groups, one can choose a basis for $K(s)^*(BG)$.  Lemma \ref{lem:free4module} ii) suggests the following

\medskip

\begin{lemma}
\label{basis-39-40}
A basis for $K^*(s)(BG)$, $G=G_{39},G_{40}$ is

$\{x_2^iy_2^j | i, j<2^{2s-1}\}$;

 $\{ax_2^iy_2^j | i<2^s, j<2^{s-1}\}$;

 $\{bx_2^iy_2^j | i<2^{2s-1},j<2^{s-1}\}$;

 $\{Tx_2^iy_2^j | i<2^{2s-1}, j<2^{s-1}(2^s-1)\}$;

 $\{c^ix_2^jy_2^k,\ c^iax_2^jy_2^k | 0<i<2^s,j<2^s,k<2^{s-1}\}$.
\end{lemma}
\begin{proof}
One can work modulo $c$ and check that first four lines give a basis for $K^*(BG)/ker \rho^*$, and then the last line forms a basis for $ker \rho^*,$ where $\rho:BH\rightarrow BG$.

Choose the lexicographic monomial ordering (lp) corresponding to the variables  $(a,T,b,y_2,x_2,c)$ in that order. Then the first four lines constitute a Gr\"{o}bner basis of
$K(s)^*(BG)/ker \rho^*$. The last line, a Gr\"{o}bner basis of $ker \rho^*$, is the union of
$\{c^ix_2^jy_2^k,\, c^iax_2^jy_2^k, | 0<i<2^s-1\}$, a basis of $ker \rho^* \bigcap K(s)^*(BG)/Im Tr^*$,
and $\{c^{2^s-1}x_2^jy_2^k,\,\,c^{2^s-1}ax_2^jy_2^k\} $, a basis of $Tr^*(\mathcal{T})$, the image of the trivial module $\mathcal{T}$ of Lemma \ref{lem:free4module} iii) under the transfer homomorphism. For the last sentence recall $Tr^*(1)=v_sc^{2^s-1}.$

Let us give the proof in the following steps and in this way explain the range restrictions for indices.

Step 1. Any monomial of $ker\rho*$ is decomposable into a sum of elements from the last line of Lemma \ref{basis-39-40}.

$cb:$ (as $cb^2=c^2b$ the decomposition of $cb^i$ will follow).
\newline
Multiply the  decomposition of $x_1$ by $c$ and take into account the relation
$c(c+x_1+\sum c^{2s-2^i}x_2^{s-1})=0$  (note also $bc^2=b^2c$ implies $c(bc)^{2^{s-1}}=0$).
This gives the decomposition of $cb$ into the $c^ix_2^j,$ $0< i<2^s, j \leq {2^{s-1}}$. Namely
\begin{equation}
\label{eq:cb}
cb=c(x_2^{2^{s-1}}+\sum_{i=1}^{s-1}c^{2^s-2^i}x_2^{2^{i-1}}).
\end{equation}

\medskip

$cx_2^{2^s}$: Multiply the decomposition of $x_2^{2^s}$  by $c$ .  As $c^{2^s}=0,$  $a^2c=ac^2$ and $b^2c=bc^2$ we have
$cx_2^{2^s}=a^2c+b^2c+ac^2=b^2c=bc^2$. Then by \eqref{eq:cb} we get

\begin{equation}
\label{eq:cx^2^s}
cx_2^{2^s}=c^2(x_2^{2^{s-1}}+\sum_{i=2}^{s-1}c^{2^s-2^i}x_2^{2^{i-1}}).
\end{equation}

$cy_2^{2^{s-1}}$: One has
\begin{equation}
\label{eq:cy^2^{s-1}}
cy_2^{2^{s-1}}=c\sum_{i=1}^{s-1}c^{2^s-2^i}y_2^{2^{i-1}}+
\begin{cases}
0,& \text{$G=G_{39}$}\\
c^2,& \text{$G=G_{40}$}
\end{cases}
\end{equation}

For this multiply the decomposition of $y_1$ of Theorem \ref{thm:G} by $c$ and apply the relation $c(c+y_1+v_s\sum_{i=1}^{s-1}c^{2^s-2^i}y_2^{2^{i-1}})=0.$

Now as $ca^2=c^2a$ and $cT=0$ we have proper decomposition for any monomial having factor $c$.

Note one has
\begin{equation}
\label{eq:x^2s-1,y^2s-1}
x_2^{2^{2s-1}}= ac^{2^s-1}, \,\,\,y_2^{2^{2s-1}}= c^{2^s-1}x_2^{2^{s-1}}.
\end{equation}

This explains the range restrictions for the first line of the basis of Lemma \ref{basis-39-40}. For this we need the decompositions of $x_2^{2^s}$ and $y_2^{2^s}$ of Theorem \ref{thm:G}.
$$x_2^{2^{2s-1}}=(ac)^{2^{s-1}}=ac^{2^s-1}.$$
Similarly
$$y_2^{2^{2s-1}}=(bc)^{2^{s-1}}=c^{2^s-1}b=c^{2^s-1}x_2^{2^{s-1}}.$$
For the last two equalities apply $cb^2=c^2b$ and \eqref{eq:cb}.

$bx_2^{2^{2s-1}}:$ Multiply the first equation of \eqref{eq:x^2s-1,y^2s-1} by $b$ and apply \eqref{eq:cb}. This gives

\begin{equation}
\label{eq:bx^2s-1}
bx_2^{2^{2s-1}}=ac^{2^s-1}x_2^{2^{s-1}}.
\end{equation}

\medskip

Step 2.

$b^2:$ Rewrite the decomposition of $y_2^{2^s}$ of Theorem \ref{thm:G} and apply \eqref{eq:cb} to get the proper decomposition
\begin{equation}
\label{eq:b2}
b^2=
\begin{cases}
y_2^{2^s}+c(x_2^{2^{s-1}}+\sum_{i=1}^{s-1}c^{2^s-2^i}x_2^{2^{i-1}}), & \text{$G=G_{39}$} \\
y_2^{2^s}+c^2+c(x_2^{2^{s-1}}+\sum_{i=1}^{s-1}c^{2^s-2^i}x_2^{2^{i-1}}), & \text{$G=G_{40}.$}
\end{cases}
\end{equation}

\medskip

$a^2:$  By Theorem \ref{thm:G} $a^2=x_2^{2^s}+b^2+ac+abc^{2^s-1}$. By \eqref{eq:cb} $abc^{2^s-1}=ac^{2^s-1}x_2^{2^{s-1}}$. Taking into account \eqref{eq:b2} we get the proper decomposition

\begin{equation}
\label{eq:a^2}
a^2=x_2^{2^s}+ac+ac^{2^s-1}x_2^{2^{s-1}}+
\begin{cases}
y_2^{2^s}+c(x_2^{2^{s-1}}+\sum_{i=1}^{s-1}c^{2^s-2^i}x_2^{2^{i-1}}), & \text{$G=G_{39}$} \\
y_2^{2^s}+c^2+c(x_2^{2^{s-1}}+\sum_{i=1}^{s-1}c^{2^s-2^i}x_2^{2^{i-1}}), & \text{$G=G_{40}.$}
\end{cases}
\end{equation}

In the following we will work modulo $c$ in the ring $R$ with lexicographic ordering determined by variables $(y_1,x_1,a,T,b,y_2,x_2,c)$ in that order and give the decompositions in the above Gr\"{o}bner basis of $K(s)^*(BG)/ker \rho^*$.

\medskip

$by_2^{2^{s-1}}:$  Clearly we need the ideal $I_1$, generated by the following relations of Theorem \ref{thm:G}: $b^{2^s},$ $b(b+y_1+\sum_{i=1}^{s-1}b^{2^s-2^i}y_2^{2^{i-1}})$, the decomposition of $y_1$ and \eqref{eq:b2} modulo $c$. Then in the quotient ring $R/I_1$
$by_2^{2^{s-1}}$ is decomposable into the elements $\{y_2^{j}\}$.

\medskip

$ab :$ Let $I_2$ be the ideal generated by the relations of $I_1$, \eqref{eq:a^2},
$a(a+x_1+\sum_{i=1}^{s-1}a^{2^s-2^i}x_2^{2^{i-1}})$ and decomposition of $x_1$ multiplied by $a$.  This gives the decomposition of $ab$ in $R/I_2$ into the elements of the first three lines of our basis.
\medskip

$bT:$ Let $I_3$ be generated by the relations of $I_2$ and $T(b+y_1+v_s\sum_{i=1}^{s-1}b^{2^s-2^i}y_2^{2^{i-1}})+b^{2^s-1}y_2(c+x_1).$
Then $bT$ is decomposable in the quotient ring $R/I_3$ into the elements of the first and fourth lines of our basis

\medskip

$T^2:$ Consider the ideal

$I_4=(I_3,T^2+Tx_1y_1+x_2y_1(c+y_1+\sum_{i=1}^{s-1}c^{2^s-2^i}y_2^{2^{i-1}})+x_1y_2(c+x_1+\sum_{i=1}^{s-1}c^{2^s-2^i}x_2^{2^{i-1}})).$
\newline
This gives a decomposition of $T^2$ which is not yet the proper decomposition, one has terms $Ty^i$, $i\geq (2^s-1)2^{s-1}.$ We will need the decomposition of $Ty_2^{2^{2s-1}-2^{s-1}}$ below to get the proper decomposition of $T^2$ into elements of the first, third and fourth lines of our basis.

\medskip

$aT:$ We need the ideal

$I_5=(I_4,T(a+x_1+\sum_{i=1}^{s-1}a^{2^s-2^i}x_2^{2^{i-1}})+a^{2^s-1}x_2(c+y_1)).$
\newline
This gives the proper decomposition of $aT$ into the first and fourth lines of our basis.

\medskip

$ay_2^{2^{s-1}}:$ Take the ideal $I_6=I_5$+the decomposition of $y_1$ multiplied by $a$  and the relation
$(c+y_1+\sum_{i=1}^{s-1}c^{2^s-2^i}y_2^{i-1})( a+x_1+\sum_{i=1}^{s-1}a^{2^s-2^i}a_2^{i-1})+a^{2^s-1}T.$
\newline
Again this gives a decomposition with terms $Ty^i$, $i\geq (2^s-1)2^{s-1}$. To get the proper decomposition into the first, third and fourth lines of our theorem we need the following decomposition of

\medskip

$Ty_2^{2^{2s}-2^{s-1}}:$ Finally put

$I_7=(I_6,(c+x_1+\sum_{i=1}^{s-1}c^{2^s-2^i}x_2^{2^{i-1}})(b+y_1+\sum_{i=1}^{s-1}b^{2^s-2^i}y_2^{2^{i-1}})+b^{2^s-1}T).$
\newline
By the above decomposition of $bT$ one has $b^{2^s-1}T=Ty_2^{2^{2s}-2^{s-1}}$ for the last summand above. This is what we need for the decomposition of $Ty_2^{2^{2s}-2^{s-1}}$ into the third and fourth lines of our basis.

Finally note that the decomposition of $ax_2^{2^s}$ already follows from the decomposition of  $x_2^{2^s}$ of Theorem \ref{thm:G} and decompositions of
$ab$ and $a^2$.

\end{proof}

Similarly Lemma \ref{lem:H38free4module} ii) suggests the following

\begin{lemma}
\label{basis38-41}

A basis for $K(s)^*(BG)$, $G=G_{38},G_{41}$ is

 $\{x_1^iy_1^jx_2^ky_2^l | i,j<2^s, k,l<2^{s-1}\}$;

 $\{abx_2^ky_2^l | k,l<2^{s-1}\}$;

 $\{y_1^iax_2^ky_2^l$, $x_1^ibx_2^ky_2^l | i<2^s, k,l<2^{s-1}\}$;

 $\{Tx_1^iy_1^jx_2^ky_2^l |i,j<2^s-1, k,l<2^{s-1}\}$;

 $\{c^ix_2^ky_2^l,\ c^iax_2^ky_2^l,\ c^ibx_2^ky_2^l,\ c^iabx_2^ky_2^l | 0<i<2^s, k,l<2^{s-1}\}$.
\end{lemma}

\medskip

Let us give a sketch of the proof. Choose the lexicographic ordering corresponding to $(T,a,b,y_2,x_2,c)$ in that order. This eliminates $a$ and $b$ (by decompositions of $x_1$ and $y_1$ of Theorem \ref{thm:G}). Then again we have to apply the relations of Theorem \ref{thm:G} and extract the following Gr\"{o}bner basis of $K(s)^*(BG)/ker \rho^*$:

\medskip

$\{x_1^iy_1^jx_2^ky_2^l | i,j<2^s, k,l<2^{s-1}\}$;

$\{x_2^{2^{s-1}+2^k}y_2^{2^{s-1}+2^l} | k,l<2^{s-1}\}$;

$\{x_2^{2^{s-1}+2^k}y_1^iy_2^l$, $x_1^i x_2^ky_2^{2^{s-1}+2^l}  | i<2^s, k,l<2^{s-1}\}$;

$\{Tx_1^iy_1^jx_2^ky_2^l |i,j<2^s-1, k,l<2^{s-1}\}$;

\medskip

(here $a$ is replaced by $x_2^{2^{s-1}}$ and $b$ by $y_2^{2^{s-1}}$ in the first four lines of Lemma \ref{basis38-41})
and a Gr\"{o}bner basis of $ker \rho^*$
$$\{c^ix_2^ky_2^l,\,\,c^ix_2^{2^{s-1}+2^k}y_2^l,\,\,c^ix_2^ky_2^{2^{s-1}+2^l},\,\,c^i x_2^{2^{s-1}+2^k}y_2^{2^{s-1}+2^l} | 0<i<2^s, k,l<2^{s-1}\}$$
corresponding to the last line of Lemma \ref{basis38-41}.

\section{Remarks}

The families of non-abelian $p$-groups whose Morava $K$-theory is known to be good in the sense of Hopkins-Kuhn-Ravenel is listed in \cite{SCH1}. In particular, if $G$ belongs to any of the following families of $p$-groups, then $K(n)^{odd}(BG) = 0.$

\medskip

(a) wreath products of the form $H \wr C_p$ with $H$ good \cite{HKR}, \cite{HU};

(b) metacyclic $p$-groups \cite{TY2};

(c) minimal non-abelian $p$-groups, i.e., groups all of whose maximal subgroups
are abelian \cite{Y1};

(d) groups of $p$-rank 2 \cite{Y2};

(e) elementary abelian by cyclic groups, i.e., the extensions $V \to G \to C$ with $V$
elementary abelian and $C$ cyclic \cite{Y3}, \cite{K};

(f) central product of the form $H\circ C_{p^m}$ with $H$ good \cite{SCH1}.

(g) $H$ is a normal subgroup in $G$ of index $p$, $H$ is good and the integral Morava $K$-theory $\tilde{K}(s)(BH)$ is a permutation module for the action of $G/H$ \cite{K}.

\medskip

For these families the ring structure of $K(s)^*(BG)$ is either studied in the works mentioned above or can be read off from previously performed computations modulo some definite indeterminacy. Namely, Yagita and Tezuka determined the multiplicative structure modulo the transfer formula \eqref{eq:tr}. On the other hand, our main aim here is to check (at least for the groups with maximal abelian subgroup of index 2) whether the transfer formula is sufficient to get the ring structure in combination with the methods of characteristic classes and transfer (double coset formula, etc.) The papers \cite{B, B1, B2, BV} treat the same problem.

Schuster suggested an alternative way to obtain explicit relations by choosing some artificial generators in the spectral sequence, not equal to Chern classes \cite{SCH1, SCHY}.

There are 51 groups of order 32. The first 7 groups are abelian and the next 8 have an abelian factor, hence the task of computing the ring structure is reduced to the smaller nonabelian groups. We refer the reader to \cite{SCH3, SCH4} for some details. In this paper, we carry out the complete
details for the groups $G$ in the Hall-Senior list with the numbers $39, \cdots ,41$.

\section*{Acknowledgements}
The authors are very grateful to the referee for exceptionally thorough analysis of the paper and numerous important suggestions which have been very useful for improving the paper.

\bibliographystyle{amsplain}

\end{document}